\documentclass{svproc}

\usepackage{graphicx}
\usepackage{amsfonts,amsmath}

\usepackage{hyperref}
\usepackage{pgfplots} 
\pgfplotsset{compat=newest}
\usetikzlibrary{plotmarks}
\usetikzlibrary{arrows.meta}
\usepgfplotslibrary{patchplots}
\usepackage{grffile}

\newcommand{\ee}{\text{e}}

\newcommand{\brho}{\boldsymbol \rho}
\newcommand{\bpsi}{\boldsymbol \psi}

\newcommand{\RR}{\mathbb{R}}
\newcommand{\CC}{\mathbb{C}}

\begin{document}
\mainmatter
\title{Exponential integrators for mean-field selective optimal control
  problems}
\author{Giacomo Albi\inst{1} \and%
  Marco Caliari\inst{2} \and%
  Elisa Calzola\inst{3} \and%
  Fabio Cassini\inst{4}}
\institute{University of Verona, \email{giacomo.albi@univr.it}\and
  University of Verona, \email{marco.caliari@univr.it}\and
  University of Verona, \email{elisa.calzola@univr.it}\and
University of Trento, \email{fabio.cassini@unitn.it}}
\maketitle

\begin{abstract}
  In this paper we consider mean-field optimal control problems
  with selective action of the control, where the constraint is a
  continuity equation involving a non-local term and diffusion.
  First order optimality conditions are formally derived in a general
  framework, accounting for boundary conditions.
  Hence, the optimality system is used to construct a reduced gradient method,
  where we introduce a novel algorithm for the numerical realization of the
  forward and the backward equations, based on exponential integrators.
  We illustrate extensive numerical experiments on different control problems
  for collective motion in the context of opinion formation and pedestrian
  dynamics.
\end{abstract}
  \section{Introduction}
  The study of collective motion of interacting agents systems is of paramount
  importance to understand the formation of coherent global behaviors
  at various scales, with applications to the study of biological, social,
  and economic phenomena. In recent years, there has been a surge of
  literature on the collective behavior of multi-agent systems,
  covering a wide range of topics such as cell aggregation and motility,
  coordinated animal motion \cite{cucker2007emergent,d2006self}, opinion
  formation \cite{MR2887663,motsch2014heterophilious,MR2247927},
  coordinated human
  behavior \cite{MR3308728,dyer2009leadership,Piccoli_2009}, and
  cooperative robots \cite{CKPP19,Meurer,KPAsurvey15,MR3157726}. These
  fields are vast and constantly evolving,
  we refer to the following surveys~\cite{bellomo20review,MR3119732,MR2580958}
  that provide a comprehensive overview of recent developments.  
Modeling such complex and diverse systems poses a significant challenge,
since in general there are no first-principles as, for instance,
in classical physics, or statistical mechanics.  Nevertheless,
the dynamics of the individuals have been successfully described by systems
of Ordinary Differential Equations (ODEs) from Newton's laws designing
basic interaction rules, such as attraction, repulsion and alignments, or,
alternatively, by considering an evolutive game where the dynamics is
driven by the simultaneous optimization of costs by $N$ players such as in
References~\cite{huang2006large,MFG}.
In this context, of paramount importance for several applications is the
design of centralized policies able to optimally enforce a desired state of
the agents, see for
instance References~\cite{MR3542027,AHP15,caponigro2015sparse}.

In this paper, we consider a constrained setting, where interacting individuals
are influenced by a centralized control with selective action, i.e.,
\begin{equation}\label{eq:min_d}
 \mathrm{d}x_i = \left( \frac{1}{N} \sum_{j =1}^N p(x_i,x_j)(x_j-x_i) +s(t,x_i,\rho^N)u_i \right) \mathrm{d}t+ \sigma \mathrm{d}W^t_i,
\end{equation}
with initial data $x^0=[x_1^0,\ldots,x_N^0]$.
Here each agent $x_i\in\Omega\subseteq \mathbb{R}^d$, for $i = 1, \ldots, N$,
accounts for pairwise interactions weighted by the function
$p(\cdot,\cdot)$, and for disturbances modelled with a Brownian motion.
The action of the control $u = [u_1, \ldots, u_N]$ is weighted by a selective
function $s(t,x_i,\rho^N)$, with $\rho^N(x)$ the empirical measures associated
to the interacting agent system, i.e.,
$\rho^N(t,x) = N^{-1}\sum_{i=1}^N \delta({x_i(t)}-x)$.
Then, the optimal control $u^*$ is obtained in the space of admissible
controls $U$, by minimizing the cost functional
\begin{equation}\label{eq:constr_d}
  J(u;x^0) =
  \mathbb{E} \left[\int_0^T \frac{1}{2N} \sum_{i =1}^N\ell(t,x_i,\rho^N)+ \gamma |u_i|^2\right],
\end{equation}
where $\ell(t,x_i,f^N)$ is a running cost to be designed by the controller,
with a quadratic penalization of the control for $\gamma\geq 0$.

For a large number of agents, we can write the mean-field optimal control
problem corresponding to the finite dimensional optimal control
problem~\eqref{eq:min_d}--\eqref{eq:constr_d}  as follows
(see References~\cite{albi2022mean,FPR14,MR3264236})
\begin{subequations}\label{eq:model0}
\begin{equation}\label{eq:min0}
 \min_{u\in U}\frac{1}{2}\int_0^T\int_\Omega\left(\ell(t,x,\rho)+\gamma|u|^2\right)\rho dxdt,     
\end{equation}
where $\rho$ is the density function satisfying the Partial
Differential Equation (PDE)
\begin{equation}\label{eq:constr0}
\left\{
  \begin{aligned}
  &\partial_t\rho+\nabla\cdot\left(\left(\mathcal{P}(\rho)+
    s(t,x,\rho)u\right)\rho\right)-\frac{\sigma^2}{2}\Delta\rho = 0,\\
    &\rho(0,x) = \rho_0(x).
 \end{aligned}\right.
\end{equation}
\end{subequations}
Here the non-local interactions among agents are described by
the integral term
\begin{equation}\label{eq:kern}
\mathcal{P}(\rho)(t,x) = \int_\Omega p(x,y)(y-x)\rho(t,y)dy
\end{equation}
and $\rho_0(x)$ is the initial distribution of the agents.
Differently from mean-field games \cite{achdou2010mean,cannarsa2021mean,MFG},
in this context the goal is to compute a mean-field optimal strategy capable
of driving the population
density to a specific target, avoiding the curse of dimensionality induced by
the large scale non-linear system of $N$ agents. However, the numerical
solution of the PDE-constrained optimization
problem~\eqref{eq:min0}--\eqref{eq:constr0} requires careful treatment
\cite{borzi2011computational}. To this end, we follow a reduced gradient method,
where the first order optimality system is solved iteratively for the
realization of the control, as in
References~\cite{aduamoah2022pseudospectral,ACFK17,bailo2018optimal}.
Major challenges arise from the presence of the stiff diffusive and
transport operators, and from the stability and storage requirements
originated by the choice of the numerical solvers. For these kinds of problems,
explicit time marching schemes usually require several time steps due to the
lack of favorable stability properties, while implicit ones need possibly
expensive solutions of (non)linear systems~\cite{albi2019linear,hager2000runge,herty2013implicit}.
A prominent and effective alternative way to numerically integrate stiff
equations in time is to employ explicit \textit{exponential integrators},
see~Reference~\cite{HO10} for a seminal review. After semidiscretization
in space, these schemes require to approximate the action of exponential
and of exponential-like matrix functions.

The paper is structured as follows. In Section~\ref{sec:continuous_model} we
present a model of interest which generalizes
the one in formulas~\eqref{eq:model0},
and we derive the formal optimality conditions using the associated Lagrangian
function, obtaining a system of coupled PDEs.
The first one is  forward in time for the density function, while
the second is backward in time for the adjoint variable.
We numerically couple these equations using the steepest descent algorithm.
In Section~\ref{numerical_integrators} we present the semidiscretization
in space of the forward and of the backward PDEs, together with the numerical
solution of the arising systems of ODEs using a pair of
exponential integrators. For convenience of the reader, we also present
there the derivation of the schemes and a brief discussion on common techniques
to compute the involved matrix functions.
Section \ref{numerical_experiments} is devoted to some numerical
validations and simulations
in opinion formation (Sznajd, Hegselmann--Krause,
and mass transfer) and pedestrian (see~Reference~\cite{BdFMW13}) models.
We finally draw some conclusions in Section~\ref{sec:conclusions}.

  \section{Mean-field selective optimal control problem}\label{sec:continuous_model}
   We consider the mean-field optimal control problem \cite{ACFK17,BdFMW13,MR3264236} defined by the functional minimization
\begin{subequations}    \label{eq:model}
   \begin{equation}\label{eq:min}
     \min_{u}\mathcal{J}(u;\rho_0),
     \end{equation}
 where $\rho = \rho(t,x)$ is a probability density of agents satisfying
   \begin{equation}\label{eq:constr}
    \left\{
    \begin{aligned}
      &\partial_t\rho+\nabla\cdot\left[\left(\mathcal{P}(\rho)+
        s(t,x,\rho)u\right)\rho\right]-\frac{\sigma^2}{2}\Delta\rho = 0, \\
    &\rho(0,x) = \rho_0(x),\\
      &\left(\left(\mathcal{P}(\rho)+
      s(t,x,\rho)u\right)\rho -\frac{\sigma^2}{2}\nabla\rho\right)\cdot \vec{n} =
    \left\{
    \begin{aligned}
    &\beta\rho & \text{on } \Gamma_\mathrm{F},\\
    & 0 & \text{on } \Gamma_\mathrm{Z}.
    \end{aligned}
    \right.
    \end{aligned}
    \right.
    \end{equation}
\end{subequations}
and defined for
 each $(t,x)\in[0,T]\times\Omega$.
 The  evolution of the density is driven by the non-local operator
 $\mathcal{P}(\rho)(t,x)$, as in equation~\eqref{eq:kern}, and by
 the control $u=u(t,x)$ 
 weighted by the selective function $s(t,x,\rho)$.
Here, we denoted by $\Gamma_\mathrm{F}$ the subset of the boundary in which
    there is a flux different from zero  ($\beta\ne 0$) and by
    $\Gamma_\mathrm{Z}$ the part of $\partial \Omega$ with zero-flux
    boundary conditions. These two subsets are such that
    $\Gamma_\mathrm{F} \cup \Gamma_\mathrm{Z} = \partial \Omega$ and
    $\Gamma_\mathrm{F} \cap \Gamma_\mathrm{Z} = \emptyset$, and $\vec{n}$ is
    the outward normal vector to the boundary with norm equal to one.
Finally, the functional in formula~\eqref{eq:min} is given by
    \begin{equation*}
      \mathcal{J}(u;\rho_0)=\frac{1}{2}\int_0^T\int_\Omega\left(e(t,x,\rho)+\gamma|u|^2\rho\right) dxdt + \frac{1}{2}\int_\Omega c(T,x,\rho(T,x))dx
    \end{equation*}
    for a general running cost $e(t,x,\rho)$ and a terminal cost $c(T,x,\rho(T,x))$.
\subsection{First order optimality conditions}
We can derive the first order optimality conditions on a formal level using
a Lagrangian approach. For a rigorous treatment we refer
to References~\cite{ACFK17,burger2021mean}. We define the Lagrangian function
with adjoint variable $\psi$ as
    \begin{equation}\label{eq:lag}
      \begin{aligned}
      \mathcal{L}(u,\rho,\psi) = &\frac{1}{2}\int_0^T\int_\Omega\left(e(t,x,\rho)+\gamma|u|^2\rho\right) dxdt + \frac{1}{2}\int_\Omega c(T,x,\rho(T,x))dx \\
      &- \int_0^T\int_\Omega\psi\left(\partial_t\rho+\nabla\cdot\left[\left(\mathcal{P}(\rho)+
        s(t,x,\rho)u\right)\rho\right]-\frac{\sigma^2}{2}\Delta\rho\right)dxdt.
    \end{aligned}
    \end{equation}
The optimal solution $(u^*,\rho^*,\psi^*)$ can be found by equating to zero the
partial Fr\'echet derivatives of the Lagrangian function, i.e.,
by solving the following system
    \begin{equation}\label{eq:optimcond}
    \left\{
    \begin{aligned}
       D_u\mathcal{L}(u,\rho,\psi)=0, \\
       D_\psi\mathcal{L}(u,\rho,\psi)=0, \\
       D_\rho\mathcal{L}(u,\rho,\psi)=0.
    \end{aligned}
    \right.
    \end{equation}
    Before computing the partial derivatives in system~\eqref{eq:optimcond},
    we integrate by parts the last term appearing in the Lagrangian
    function~\eqref{eq:lag} and we get
    \begin{equation*}
\begin{split}
      \mathcal{L}(u,\rho, \psi)=&
 \frac{1}{2} \int_0^T\int_\Omega  \left(e (t,x,\rho) + \gamma |u|^2 \rho\right) d x d t
+\frac{1}{2}\int_\Omega c(T,x,\rho(T,x))dx\\
& + \int_0^T \int_\Omega\rho\left(\partial_t \psi +\frac{\sigma^2}{2}\Delta \psi+\left(\mathcal{P}(\rho)+ s(t,x,\rho)u \right) \cdot \nabla \psi\right) dxd t\\
& -\int_0^T \int_{\Gamma_\mathrm{F}}\rho\left( \frac{\sigma^2}{2}\nabla \psi \cdot \vec{n}  +
\beta\psi \right) d b d t\\
&-\int_\Omega (\psi(T,x)\rho(T,x)-\psi(0,x)\rho(0,x))dx,
\end{split}
    \end{equation*}
where we used the value of the boundary conditions appearing in equation~\eqref{eq:constr}.
 Performing then the computations of the
 partial derivatives we obtain the gradient direction
 for the control variable
 $u$
    \begin{equation}
      D_u\mathcal{L}(u,\rho,\psi)=\gamma u + s(t,x,\rho)\nabla\psi,
      \label{eq:contr}
    \end{equation}
the forward PDE for the density function $\rho$
    \begin{equation}
    \left\{
    \begin{aligned}
      &\partial_t\rho+\nabla\cdot\left[\left(\mathcal{P}(\rho)+
        s(t,x,\rho)u\right)\rho\right]-\frac{\sigma^2}{2}\Delta\rho = 0, \\
    &\rho(0,x) = \rho_0(x),\\
      &\left(\left(\mathcal{P}(\rho)+
      s(t,x,\rho)u\right)\rho -\frac{\sigma^2}{2}\nabla\rho\right)\cdot \vec{n} =
    \left\{
    \begin{aligned}
    &\beta\rho & \text{on } \Gamma_\mathrm{F},\\
    & 0 & \text{on } \Gamma_\mathrm{Z},
    \end{aligned}
    \right.
    \end{aligned}
    \right.
    \label{eq:fwd}
    \end{equation}
and the backward PDE for the adjoint variable $\psi$
    \begin{equation}
    \left\{
    \begin{aligned}
    & -\partial_t\psi = \frac{\sigma^2}{2}\Delta\psi
    +\left(\mathcal{P}(\rho)+(s(t,x,\rho)+\rho D_\rho s(t,x,\rho))u\right)\cdot\nabla\psi+\\
    &\quad\quad\quad+\mathcal{Q}(\rho,\psi) + \frac{1}{2}(D_\rho e(t,x,\rho) + \gamma |u|^2 ), \\
    &\psi(T,x) = \psi_T(x), \\
    &\frac{\sigma^2}{2}\nabla\psi \cdot \vec{n} = 
    \left\{
    \begin{aligned}
    &-\beta\psi & \text{on } \Gamma_\mathrm{F},\\
    & 0 & \text{on } \Gamma_\mathrm{Z},
    \end{aligned}
    \right.
    \end{aligned}
    \right.
    \label{eq:bwd}
    \end{equation}
    where
    \begin{equation*}
      \mathcal{Q}(\rho,\psi)(t,x) = \int_\Omega p(y,x)(x-y)\cdot\nabla\psi(t,y)\rho(t,y)dy
    \end{equation*}
    and $\psi_T(x) = \frac{1}{2}D_\rho c(T,x,\rho(T,x))$.
    Now, in order to solve model~\eqref{eq:model}, we employ a steepest descent approach
    (see References~\cite{ACFK17,bailo2018optimal}).
    Starting with an initial control $u^0$, at each
    iteration $\ell$ we insert $u^\ell$ into
    the forward equation~\eqref{eq:fwd} and solve it for
    $\rho=\rho^{\ell+1}$. We then
    insert $u^\ell$ and $\rho^{\ell+1}$ into the backward
    equation~\eqref{eq:bwd} and solve
    it for $\psi=\psi^{\ell+1}$. We finally update the control
    by using the gradient
    direction~\eqref{eq:contr}, i.e.,
    \begin{equation*}
      u^{\ell+1} = u^{\ell}-\lambda^\ell(\gamma u^\ell + s(t,x,\rho^{\ell+1})\nabla\psi^{\ell+1})
    \end{equation*}
    and get $u^{\ell+1}$.
    We proceed iterating until
    $\mathcal{J}(u^{\ell+1})$ has stabilized within a given
    tolerance.
    For the numerical solution of equations~\eqref{eq:fwd} and \eqref{eq:bwd}
    we use the method of lines: we first discretize in space
    and then use appropriate integrators for the obtained systems of ODEs.
\section{Numerical integrators for the semidiscretized equations}\label{numerical_integrators}
In this section, we explain how to solve the forward and the backward
PDEs in the steepest descent algorithm.
By observing that both are semilinear parabolic equations, the idea
is to use numerical schemes tailored for this type of problems. A prominent
way is to apply explicit exponential integrators~\cite{HO10} to the systems
of ODEs arising from the semidiscretization in space of the PDEs.
By construction, these schemes solve exactly linear ODEs systems with constant
coefficients,
they allow for time steps usually much larger than those
required by classical explicit methods (i.e., typically they do not suffer
from a CFL restriction), and do not require the solution
of (non)linear systems as implicit methods do.
On the other hand, this class of integrators requires the computation of
the action of exponential-like matrix functions for which different efficient
techniques have been developed in recent years.
\subsection{Forward PDE}
For sake of clarity, and since we will present later on one-dimensional
numerical examples, we consider $\Omega=[a,b]$ and
we rewrite the forward PDE~\eqref{eq:fwd}
\begin{equation*}
\left\{
\begin{aligned}
  &\partial_t\rho(t,x)=\frac{\sigma^2}{2}\partial_{xx}\rho(t,x)-
  \partial_x\left((\mathcal{P}(\rho(t,\cdot))(t,x)+
    s(t,x,\rho(t,x))u(t,x))\rho(t,x)\right),\\
&\rho(0,x) = \rho_0(x),\\
    &\left((\mathcal{P}(\rho(t,\cdot))(t,x)+s(t,x,\rho(t,x))u(t,x))\rho(t,x)
    -\frac{\sigma^2}{2}\partial_x\rho(t,x)\right)\bigg|_a=
\beta_a\rho(t,a),\\
&\left((\mathcal{P}(\rho(t,\cdot))(t,x)+s(t,x,\rho(t,x))u(t,x))\rho(t,x)
-\frac{\sigma^2}{2}\partial_x\rho(t,x)\right)\bigg|_b=
\beta_b\rho(t,b),
\end{aligned}
\right.
\end{equation*}
where $\beta_a,\beta_b\in\RR$ can be selected so that it is possible to
express both zero and nonzero fluxes.
Notice that when we solve this equation we
consider $u(t,x)$ a given function. We introduce a semidiscretization in space
by finite differences on a grid of points
$x_i$, with $i=1,\ldots,n$, in such a way that
$\brho(t)=[\rho_1(t),\ldots,\rho_n(t)]^{\sf T}$ is the unknown vector
whose components
$\rho_i(t)$ approximate $\rho(t,x_i)$.
Now, by denoting $D_1$ and $D_2$
the matrices which discretize $\partial_x$ and
$\partial_{xx}$ at the grid points, respectively, and $P$ the discretization
of the linear integral operator $\mathcal{P}$ by a quadrature formula,
the linear part of the right hand side
of the equation is discretized by
\begin{equation*}
\tilde A_\mathrm{F}\brho(t)=\frac{\sigma^2}{2}D_2\brho(t),
\end{equation*}
while the nonlinear part becomes
\begin{multline*}
  \tilde {\boldsymbol g}_\mathrm{F}(t,\brho(t))=
  -(D_1P\brho(t))\brho(t)
  -(P\brho(t))(D_1\brho(t))\\
  -(D_1\boldsymbol s(t,\brho(t)))\boldsymbol u(t)\brho(t)
  -\boldsymbol s(t,\brho(t))(D_1\boldsymbol u(t))\brho(t)
  -\boldsymbol s(t,\brho(t))\boldsymbol u(t)(D_1\brho(t)).
\end{multline*}
Now, we also discretize the boundary conditions
with finite differences by using virtual
nodes, and we modify
accordingly both the linear part $\tilde A_\mathrm{F}$ and the nonlinear one
$\tilde {\boldsymbol g}_\mathrm{F}(t,\brho(t))$. The resulting nonlinear
system of ODEs is then
    \begin{equation}\label{eq:ODEFP}
      \left\{
      \begin{aligned}
        \brho'(t) &= A_\mathrm{F} \brho(t) +
        \boldsymbol g_\mathrm{F}(t,\brho(t)),\quad t\in[0,T],\\
      \brho(0)&=\brho_0.
      \end{aligned}\right.
    \end{equation}

Given a time discretization $[t_0,\ldots,t_k,\ldots,t_m]$, with $t_0=0$ and
$t_m=T$, the exact solution of system~\eqref{eq:ODEFP}
at time $t_{k+1}$ can be expressed using the
variation-of-constants formula, i.e.,
\begin{equation*}
  \brho(t_{k+1}) = \ee^{\tau_{k+1} A_\mathrm{F}}\brho(t_k) +
  \int_0^{\tau_{k+1}}\ee^{(\tau_{k+1}-s)A_\mathrm{F}}
  \boldsymbol g_\mathrm{F}(t_k+s,\brho(t_k+s))ds,
\end{equation*}
where $\tau_{k+1}=t_{k+1}-t_k$, for $k=0,\ldots,m-1$.
In order to obtain an explicit
first order numerical scheme, we denote by $\brho_k$
the approximation of $\brho(t_k)$ and approximate
the nonlinear function $\boldsymbol g_\mathrm{F}(t_k+s,\brho(t_k+s))$
with $\boldsymbol g_\mathrm{F}(t_k,\brho_k)$. Hence, we have
\begin{equation}\label{eq:expEuler}
  \begin{aligned}
    \brho(t_{k+1})\approx\brho_{k+1} &= \ee^{\tau_{k+1} A_\mathrm{F}}\brho_k +
    \int_0^{\tau_{k+1}}\ee^{(\tau_{k+1}-s)A_\mathrm{F}}
    \boldsymbol g_\mathrm{F}(t_k,\brho_k)ds \\
    &= \ee^{\tau_{k+1} A_\mathrm{F}}\brho_k +
    \left(\int_0^{\tau_{k+1}}\ee^{(\tau_{k+1}-s)A_\mathrm{F}}ds\right)
    \boldsymbol g_\mathrm{F}(t_k,\brho_k) \\
    &= \ee^{\tau_{k+1} A_\mathrm{F}}\brho_k +
    \left(\tau_{k+1}\int_0^1\ee^{\tau_{k+1}(1-\theta)A_\mathrm{F}}d\theta\right)
    \boldsymbol g_\mathrm{F}(t_k,\brho_k) \\
    &= \ee^{\tau_{k+1} A_\mathrm{F}}\brho_k
    + \tau_{k+1}\varphi_1(\tau_{k+1} A_\mathrm{F})
    \boldsymbol g_\mathrm{F}(t_k,\brho_k).
      \end{aligned}
\end{equation}
Here 
we introduced the exponential-like function
    \begin{equation*}
      \varphi_1(X) = \int_0^1\ee^{(1-\theta)X}d\theta,
    \end{equation*}
 with $X\in\CC^{n\times n}$ a generic matrix. This scheme is known as
\emph{exponential Euler},
it is a fully explicit method of first (stiff) order and it
is A-stable by construction.
Its implementation requires at each time step
the evaluation of a linear combination of type
$\ee^{\tau_{k+1} X}\boldsymbol v_k+\tau_{k+1}\varphi_1(\tau_{k+1} X)\boldsymbol w_k$,
where $\boldsymbol v_k,\boldsymbol w_k\in\CC^n$ are suitable vectors,
which we will address in Section~\ref{sec:mfe}.
\subsubsection{Selective function independent of the density}
A remarkable occurrence in the literature is the one in which the selective
function does not depend on the density, i.e.,
$s(t,x,\rho(t,x))=s(t,x)$
(see Reference~\cite{ACFK17} for the case $s(t,x)=1$, which we will also
consider in the numerical examples).
In this case, some terms in the nonlinear
part $\tilde{\boldsymbol g}_\mathrm{F}(t,\brho(t))$ can actually be incorporated
into the linear one. In fact, we obtain
\begin{equation*}
\tilde A_\mathrm{F}(t)\brho(t)=\frac{\sigma^2}{2}D_2\brho(t)
-(D_1\boldsymbol s(t))\boldsymbol u(t)\brho(t)
- \boldsymbol s(t)(D_1 \boldsymbol u (t))\boldsymbol \rho(t)
- \boldsymbol s(t)\boldsymbol u(t)(D_1 \boldsymbol \rho(t)),
\end{equation*}
while the nonlinear part is now given by
\begin{equation*}
  \tilde {\boldsymbol g}_\mathrm{F}(t,\brho(t))=
  -(D_1P\brho(t))\brho(t)
  -(P\brho(t))(D_1\brho(t)).
\end{equation*}
By modifying accordingly the quantities in order to impose the boundary
conditions, we end up with the system of ODEs 
    \begin{equation}\label{eq:ODEFP1}
      \left\{
      \begin{aligned}
        \brho'(t) &= A_\mathrm{F}(t) \brho(t) +
        \boldsymbol g_\mathrm{F}(t,\brho(t)),\quad t\in[0,T],\\
      \brho(0)&=\brho_0,
      \end{aligned}\right.
    \end{equation}
which is similar to system~\eqref{eq:ODEFP}, except for the fact that the linear
part has time dependent coefficients. Nevertheless,
at each $t_k$ we can rewrite equivalently this system as
    \begin{equation*}
      \left\{
      \begin{aligned}
        \brho'(t) &= A_\mathrm{F}(t_k) \brho(t) +
        ( A_\mathrm{F}(t)- A_\mathrm{F}(t_k))\brho(t)+
        \boldsymbol g_\mathrm{F}(t,\brho(t))\\
        & = A_\mathrm{F}(t_k) \brho(t) +
        \boldsymbol g^{k}_\mathrm{F}(t,\brho(t)),\\
      \brho(0)&=\brho_0.
      \end{aligned}\right.
    \end{equation*}
and apply the exponential Euler method. Thus,  we end up with the scheme
\begin{equation}\label{eq:expEulerMagnus1}
  \begin{split}
  \brho(t_{k+1})\approx\brho_{k+1}&=\ee^{\tau_{k+1} A_\mathrm{F}(t_{k})}\brho_{k}+
  \tau_{k+1}\varphi_1(\tau_{k+1} A_\mathrm{F}(t_{k}))
  \boldsymbol g_{\mathrm{F}}^{k}(t_{k},\brho_{k})\\
  &=\ee^{\tau_{k+1} A_\mathrm{F}(t_{k})}\brho_{k}+
  \tau_{k+1}\varphi_1(\tau_{k+1} A_\mathrm{F}(t_{k}))
  \boldsymbol g_{\mathrm{F}}(t_{k},\brho_k),
  \end{split}
\end{equation}
for $  k=0,\ldots,m-1$.
As for the general case $s(t,x,\rho(t,x))$, we obtain in this way an explicit
method of first order (which we call \emph{exponential Euler--Magnus})
that requires again a linear combination of actions of the matrix exponential
and the matrix $\varphi_1$ function.
\subsection{Backward PDE}
We rewrite the backward PDE~\eqref{eq:bwd} in the one-dimensional
case $\Omega=[a,b]$
\begin{equation*}
\left\{
\begin{aligned}
  &-\partial_t\psi(t,x) = \frac{\sigma^2}{2}\partial_{xx}\psi(t,x)
  +\mathcal{P}(\rho(t,\cdot))(t,x)\partial_x\psi(t,x)\\
  &\phantom{-\partial_t\psi(t,x) =\;\; }
  +(s(t,x,\rho(t,x))+\rho(t,x)s_\rho(t,x,\rho(t,x)))u(t,x)\partial_x\psi(t,x)\\
  &\phantom{-\partial_t\psi(t,x) =\;\; }
  +\mathcal{Q}(\rho(t,\cdot),\psi(t,\cdot))(t,x)+\frac{1}{2}\left(e_\rho(t,x,\rho(t,x)) + \gamma u^2(t,x) \right), \\
&\psi(T,x) = \psi_T(x), \\
&\frac{\sigma^2}{2}\partial_x\psi(t,x)\big|_a =
-\beta_a\psi(t,a),\\
&\frac{\sigma^2}{2}\partial_x\psi(t,x)\big|_b =
-\beta_b\psi(t,b),
\end{aligned}
\right.
\end{equation*}
where $s_\rho(t,x,\rho(t,x))=D_\rho s(t,x,\rho(t,x))$ and
$e_\rho(t,x,\rho(t,x))=D_\rho e(t,x,\rho(t,x))$.
Here we assume that $\rho(t,x)$ and $u(t,x)$ are given functions.
By applying a finite difference discretization on the same spatial grid
as above and defining
$Q$ the discretization of the linear integral  operator
$\mathcal{Q}$ we obtain
the linear part
\begin{multline*}
  \tilde A_{\mathrm{B}}(t)\bpsi(t)=
  \frac{\sigma^2}{2}D_2\bpsi(t)+
  (P\brho(t))(D_1\bpsi(t))\\
  +(\boldsymbol s(t,\brho(t))+
  \brho(t)\boldsymbol s_{\brho}(t,\brho(t)))\boldsymbol u(t)(D_1\bpsi(t))
  +Q(\brho(t)(D_1\bpsi(t)))
\end{multline*}
and the source term
\begin{equation*}
  \tilde{\boldsymbol g}_{B}(t)=\frac{1}{2}\boldsymbol e_{\brho}(t,\brho(t))+
  \gamma \boldsymbol u^2(t).
\end{equation*}
Finally, by taking into consideration boundary conditions, we end up with
the inhomogeneous time dependent coefficient linear system of ODEs
\begin{equation}\label{eq:ODEHJB}
  \left\{\begin{aligned}
  -\bpsi'(t) &= A_\mathrm{B}(t) \bpsi(t) + \boldsymbol g_\mathrm{B}(t),\quad
  t\in[0,T],\\
  \bpsi(T)&=\bpsi_T.
  \end{aligned}\right.
\end{equation}
By considering the same  time discretization $[t_0,\ldots,t_{k+1},\ldots,t_m]$
as above,
system~\eqref{eq:ODEHJB} has a similar structure to system~\eqref{eq:ODEFP1}.
Hence, taking into account that we are marching
backward in time, we apply the exponential Euler--Magnus method
and we obtain the time marching
\begin{equation}\label{eq:expEulerMagnus}
   \bpsi(t_k)\approx\bpsi_{k}=\ee^{\tau_{k+1} A_\mathrm{B}(t_{k+1})}\bpsi_{k+1}+
  \tau_{k+1}\varphi_1(\tau_{k+1} A_\mathrm{B}(t_{k+1}))
  \boldsymbol g_{\mathrm{B}}(t_{k+1}),
\end{equation}
for $  k=m-1,m-2\ldots,0$.
\subsection{Matrix functions evaluation}\label{sec:mfe}
We have introduced two exponential integrators that require, at each time
step, the evaluation of
\begin{equation}\label{eq:lincomb}
  \ee^{\tau X}\boldsymbol v+
  \tau\varphi_1(\tau X)\boldsymbol w,
\end{equation}
where $\tau>0$,
$X\in\RR^{n\times n}$, and $\boldsymbol v,\boldsymbol w\in\RR^n$.
We stress that these quantities depend in general on the current
time step, but for simplicity of notation we dropped the subscripts.
If we choose a uniform time discretization, i.e., $\tau_k=\tau$ for
$k=0,\ldots,m-1$,
in the exponential Euler scheme~\eqref{eq:expEuler}
we can compute once and for all the matrices
$\ee^{\tau A_\mathrm{F}}$ and $\varphi_1(\tau A_\mathrm{F})$
and then multiply by the corresponding vectors.
In this case, for the matrix function approximations
the most common techniques are Taylor expansions or Pad\'e
rational approximations with scaling and squaring
(see, for instance, References~\cite{AMH09,CZ19,SID19,SW09}).
This approach is computationally attractive only for matrices of moderate
size, taking into account also that the resulting matrix functions are full
even if the original ones were sparse.
When employing the exponential Euler--Magnus schemes~\eqref{eq:expEulerMagnus1}
and~\eqref{eq:expEulerMagnus},
we can still pursue this approach. However, since here
the matrices change at each time step, we need to recompute the
matrix functions every time accordingly.
It is also possible to compute 
linear combination~\eqref{eq:lincomb} by using a \emph{single}
slightly augmented matrix function evaluation.
In fact, thanks to~\cite[Proposition 2.1]{SA92},
we have that the first $n$ rows of
\begin{equation*}
  \exp\left(\tau\begin{bmatrix}
    X&\boldsymbol w\\
    0\cdots0 & 0
  \end{bmatrix}\right)\begin{bmatrix}
    \boldsymbol v\\
    1
    \end{bmatrix}
\end{equation*}
coincide with vector~\eqref{eq:lincomb}. This is an attractive choice
in a variable step size scenario, in which both the forward and the backward
equations could be solved by a single matrix function evaluation at each time
step.

When $X$ is a large sized and sparse matrix, it may be convenient to compute
directly vector~\eqref{eq:lincomb} at each time step \emph{without} explicitly
computing the matrix exponential.
State-of-the-art techniques follow this approach and are based on
Krylov methods or direct interpolation polynomial methods (see, for instance,
References~\cite{AMH11,CCZ22,GRT18,LPR19}).

\section{Numerical experiments}\label{numerical_experiments}
We present in this section several numerical examples arising from
different choices of parameters and functions
in the continuous model~\eqref{eq:model}. In particular, we consider
numerical experiments for two different classes of multi-agent systems
in opinion formation and pedestrian dynamics.
In all cases, we discretize in space with second order centered finite
differences and we employ the trapezoidal rule for the quadrature of the
integral operators.
All the numerical experiments have been performed on an
Intel\textsuperscript{\textregistered}
Core\textsuperscript{\texttrademark} i7-10750H
CPU with six physical cores and 16GB of RAM, using
\textsc{matlab} programming language. As a software, we use
MathWorks MATLAB\textsuperscript{\textregistered}
R2022a.
In order to compute the needed actions of exponential and
$\varphi_1$-function, we employ the \texttt{kiops}
function\footnote{\url{https://gitlab.com/stephane.gaudreault/kiops/-/tree/master/}},
which is based on the Krylov method
and whose underlying algorithm is thoroughly presented in
Reference~\cite{GRT18}.
This routine requires an input tolerance, which we set
sufficiently small in order not to affect the accuracy of the
temporal integration.

\subsection{Control in opinion dynamics}
In this section we consider two models for control of opinion dynamics,
namely the Sznajd and the Hegselmann--Krause (bounded confidence) ones,
similarly to References~\cite{ACFK17,hegselmann2002opinion,sznajd2000opinion}.
We set both models in
the spatial domain $\Omega=[-1,1]$, whose boundaries
represent the extremal opinions. The running cost is
$e(t,x,\rho) = |x-x_d|^2\rho(t,x)$ and the selective function $s(t,x,\rho)$ is set to
the constant~1 (hence, we use the exponential Euler--Magnus
scheme~\eqref{eq:expEulerMagnus1} for the forward equation).
For both the problems we consider
in  model~\eqref{eq:model} zero-flux boundary conditions everywhere
and null terminal cost function $c(T,x,\rho(T,x))$.
\subsubsection{Sznajd model}\label{sec:Sznajd}
In the first numerical experiment we present an example of
Sznajd model for opinion formation taken from Reference~\cite{ACFK17}.
In particular,  we consider
the interaction function $p(x,y)=x^2-1$, representing a repulsive interaction,
and
the target point in the running cost $x_d=-0.5$. Moreover,
we set the penalization parameter $\gamma=0.5$ and the diffusion coefficient
$\sigma=\sqrt{0.02}$. The initial density
function is of bimodal type
\begin{equation*}
\rho_0(x) = C(\rho_+(x+0.75;0.05,0.5) + \rho_+(x-0.5;0.15,1)),
\end{equation*}
where
\begin{equation*}
\rho_+(x;a,b) = \mathrm{max}\left\{-\left(\frac{x}{b}\right)^2+a,0\right\}
\end{equation*}
and $C$ defined so that $\int_\Omega\rho_0(x)dx=1$.

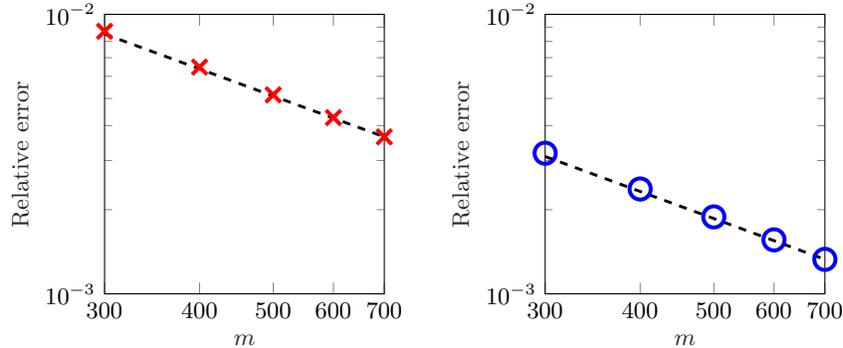
\begin{figure}[htb!]
  \centering
%
%
\begin{tikzpicture}

\begin{axis}[%
width=1.464in,
height=1.464in,
at={(0.769in,0.477in)},
scale only axis,
xmode=log,
xmin=300,
xmax=700,
xminorticks=true,
xtick={300,400,500,600,700},
xticklabels={{300},{400},{500},{600},{700}},
xlabel style={font=\color{white!15!black}},
xlabel={$m$},
ylabel style={font=\color{white!15!black}},
ylabel={Relative error},
ymode=log,
ymin=1e-3,
ymax=1e-2,
yminorticks=true,
axis background/.style={fill=white},
]
\addplot [color=red, line width=1.5pt, only marks, mark size=4pt, mark=x, mark options={solid, red}, forget plot]
  table[row sep=crcr]{%
300	8.6955e-3\\
400	6.4761e-3\\
500	5.1516e-3\\
600	4.2716e-3\\
700	3.6445e-3\\
};
\addplot [color=black, dashed, line width=1.1pt, forget plot]
  table[row sep=crcr]{%
300	8.5039e-3\\
700	3.6445e-3\\
};
\end{axis}

\begin{axis}[%
width=1.464in,
height=1.464in,
at={(3.075in,0.477in)},
scale only axis,
xmode=log,
xmin=300,
xmax=700,
xtick={300,400,500,600,700},
xticklabels={{300},{400},{500},{600},{700}},
xminorticks=true,
xlabel style={font=\color{white!15!black}},
xlabel={$m$},
ylabel style={font=\color{white!15!black}},
ylabel={Relative error},
ymode=log,
ymin=1e-3,
ymax=1e-2,
yminorticks=true,
axis background/.style={fill=white},
]
\addplot [color=blue, line width=1.5pt, only marks, mark size=4pt, mark=o, mark options={solid, blue}, forget plot]
  table[row sep=crcr]{%
300	3.1846e-3\\
400	2.3724e-3\\
500	1.8851e-3\\
600	1.5603e-3\\
700	1.3283e-3\\
};
\addplot [color=black, dashed, line width=1.1pt, forget plot]
  table[row sep=crcr]{%
300	3.0994e-3\\
700	1.3283e-3\\
};
\end{axis}

\end{tikzpicture}%
  \caption{Relative errors in infinity norm
    of $\brho(T)$ (left, $T=4$) and $\bpsi(0)$ (right),
    with respect to a reference solution, for the
    Sznajd model described in Section~\ref{sec:Sznajd} with $n=200$
    spatial discretization points and varying number of time steps $m$.
    The reference line of order 1 is also displayed.}
\label{fig:orderSznajd}
\end{figure}

First of all, we show that the expected temporal rate of convergence
of the exponential integrators is preserved also after a complete
solution of the model.
In fact, for a semidiscretization in space with $n=200$
uniform grid points, we solve
several times model~\eqref{eq:model} by the steepest
descent method described at the end of Section~\ref{sec:continuous_model}
by employing an increasing sequence
of time steps, ranging from $m=300$ to $m=700$.
Each time, after the stabilization of the functional $\mathcal{J}$,
we measure the error at final time $T=4$ for 
$\brho(t)$ and at initial time for
$\bpsi(t)$ with respect to reference solutions.
We display in Figure~\ref{fig:orderSznajd} the obtained relative errors,
which confirm the expected accuracy and rate of convergence.

\begin{figure}[htb!]
  \centering
  \includegraphics[scale=0.22]{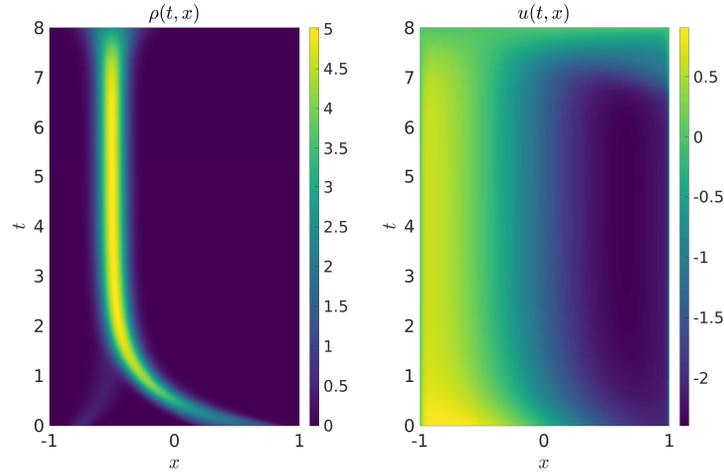}
  \caption{Evolution of the density $\rho(t,x)$ (left) and of the
    control $u(t,x)$ (right) up to final time $T=8$ for the Sznajd model
    described in Section~\ref{sec:Sznajd} with $n=1000$ spatial
    discretization points and $m=200$ time steps.}
  \label{fig:behaviorSznajd}
\end{figure}

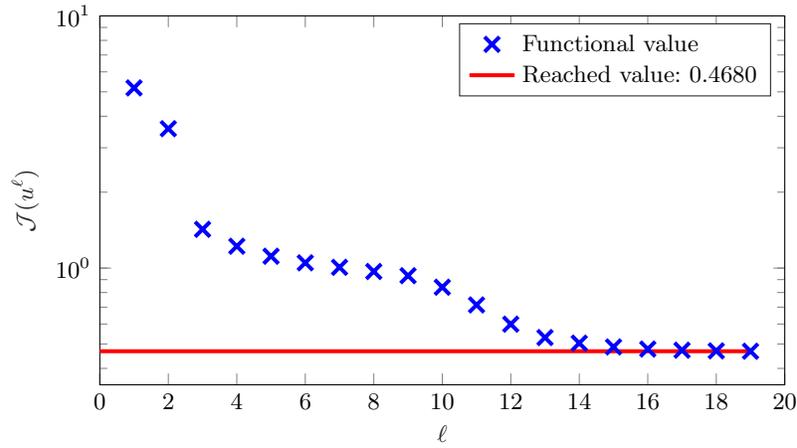
\begin{figure}[htb!]
  \centering
%
%
\begin{tikzpicture}

\begin{axis}[%
width=3.585in,
height=1.932in,
at={(0.769in,0.477in)},
scale only axis,
xmin=0,
xmax=20,
xlabel style={font=\color{white!15!black}},
xlabel={$\ell$},
ymin=0,
ymax=10,
ymode=log,
ylabel style={font=\color{white!15!black}},
ylabel={$\mathcal{J}(u^\ell)$},
axis background/.style={fill=white},
legend style={legend cell align=left, align=left, draw=white!15!black}
]
\addplot [color=blue, line width=1.5pt, only marks, mark size=4pt, mark=x, mark options={solid, blue}]
  table[row sep=crcr]{%
1	5.18314867120775\\
2	3.57096691006941\\
3	1.42492675044067\\
4	1.22098713645356\\
5	1.11418581979697\\
6	1.04919043789193\\
7	1.00716199620951\\
8	0.969412240455384\\
9	0.932860810559546\\
10	0.839626762819307\\
11	0.714441774129739\\
12	0.599134091979629\\
13	0.530313300530814\\
14	0.503871949484608\\
15	0.486652726171251\\
16	0.477379459807526\\
17	0.472457495852208\\
18	0.469640250735839\\
19	0.467979704036818\\
};
\addlegendentry{Functional value}

\addplot [color=red, line width=1.5pt]
  table[row sep=crcr]{%
0	0.467979704036818\\
19	0.467979704036818\\
};
\addlegendentry{Reached value: 0.4680}

\end{axis}

\end{tikzpicture}%
  \caption{Value of the functional  $\mathcal{J}(u^\ell)$ at the successive
iterations of the steepest descent method for the Sznajd model
described in Section~\ref{sec:Sznajd} ($n=1000$ and $m=200$).}
    \label{fig:functionalSznajd}
\end{figure}

Then, we show the behavior of the Sznajd model in opinion formation. For this
purpose we use a spatial discretization of $n=1000$ points and $m=200$ time
steps. Notice that we can employ a number of time steps small with respect
to the number of discretization points since the exponential integrators
applied to this problem do not exhibit any CFL restriction, in contrast
to explicit methods.  In Figure~\ref{fig:behaviorSznajd} we show
the evolution of the density $\rho(t,x)$ and of the control $u(t,x)$. The
results have the expected behavior of concentration of the opinions
around the target point $x_d=-0.5$ and qualitatively
match the analogous simulation available in the literature~\cite{ACFK17}.
Moreover, we show in Figure~\ref{fig:functionalSznajd} the value
of the functional $\mathcal{J}(u^\ell)$ at the successive
iterations of the steepest descent method. We observe that
the method needs 19 iterations to reach the input tolerance $2\cdot 10^{-3}$.
Finally, the overall computational time of this simulation is
about 55 seconds.

\subsubsection{Hegselmann--Krause model}\label{sec:HK}
\begin{figure}[htb!]
  \centering
  \includegraphics[scale=0.22]{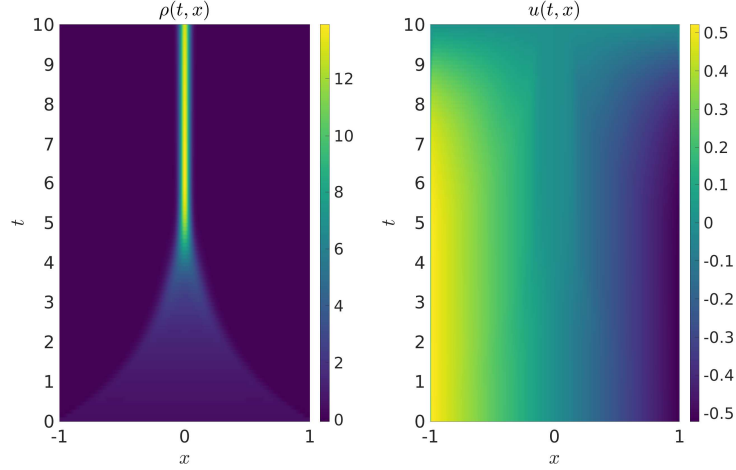}
  \caption{Evolution of the density $\rho(t,x)$ (left) and of the
    control $u(t,x)$ (right) up to final time $T=10$ for the
    Hegselmann--Krause model
    described in Section~\ref{sec:HK} with $n=1000$ spatial
    discretization points and $m=100$ time steps.}
\label{fig:behaviorHK}
\end{figure}

\begin{figure}[htb!]
  \centering
%
%
%
\begin{tikzpicture}

\begin{axis}[%
width=3.585in,
height=1.932in,
at={(0.769in,0.477in)},
scale only axis,
xmin=0,
xmax=15,
xlabel style={font=\color{white!15!black}},
xlabel={$\ell$},
ymin=0,
ymax=10,
ymode=log,
ylabel style={font=\color{white!15!black}},
ylabel={$\mathcal{J}(u^\ell)$},
axis background/.style={fill=white},
legend style={legend cell align=left, align=left, draw=white!15!black}
]
\addplot [color=blue, line width=1.5pt, only marks, mark size=4pt, mark=x, mark options={solid, blue}]
  table[row sep=crcr]{%
1	1.65549133917979\\
2	0.960355634900682\\
3	0.633548789474009\\
4	0.413235673051207\\
5	0.34476511541026\\
6	0.319068824537322\\
7	0.304041112251507\\
8	0.293988240250676\\
9	0.286811538324599\\
10	0.281465610328747\\
11	0.277356948699259\\
12	0.27412698083095\\
13	0.271484646793505\\
14	0.269349720862304\\
15	0.267671802411463\\
};
\addlegendentry{Functional value}

\addplot [color=red, line width=1.5pt]
  table[row sep=crcr]{%
0	0.267671802411463\\
15	0.267671802411463\\
};
\addlegendentry{Reached value: 0.2677}

\end{axis}

\end{tikzpicture}%
  \caption{Value of the functional $\mathcal{J}(u^\ell)$ at the successive
iterations of the steepest descent method for the Hegselmann--Krause model
described in Section~\ref{sec:HK}
($n=1000$ and $m=100$).}
\label{fig:functionalHK}
\end{figure}
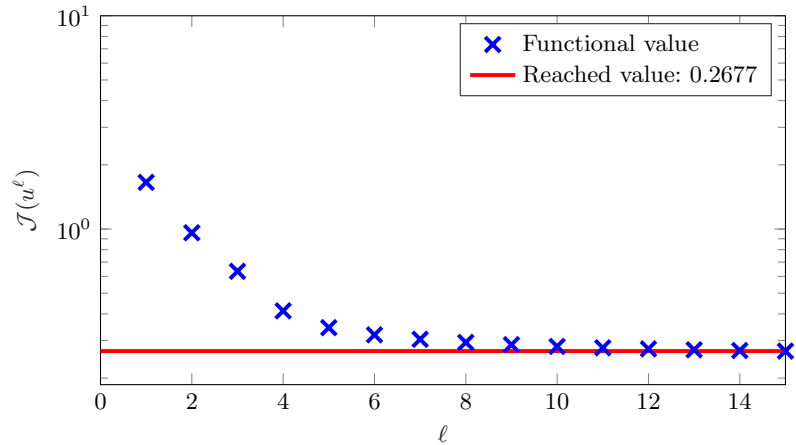

In the second numerical experiment we present an example of
Hegselmann--Krause model for opinion formation taken from
Reference~\cite{ACFK17}.
In particular,  we take
the interaction function $p(x,y)=\chi_{\{|x-y|\leq\kappa\}}(y)$,
with $\kappa=0.15$, and
the target point in the running cost $x_d=0$. Moreover,
we set the penalization parameter $\gamma=2.5$ and the diffusion coefficient
$\sigma=\sqrt{0.002}$. The initial density
function is
\begin{equation*}
\rho_0(x) = C(0.5+\epsilon(1-x^2)),
\end{equation*}
where $\epsilon=0.01$
and $C$ defined so that $\int_\Omega\rho_0(x)dx=1$.
For this model, we directly present the results
 using a spatial discretization of $n=1000$ points and $m=100$ time
 steps up to the final time $T=10$.
  In Figure~\ref{fig:behaviorHK} we display
  the evolution of the density $\rho(t,x)$ and of the control $u(t,x)$.
  Similarly to the Sznajd model, the results match both the expectations
  and the outcomes in the literature.
Then, we display in Figure~\ref{fig:functionalHK} the value
of the functional $\mathcal{J}(u^\ell)$ at the successive
iterations of the steepest descent method. We observe that
the method needs 15 iterations to reach the input tolerance $2\cdot 10^{-3}$.
Finally, this simulation takes roughly 15 seconds.

\subsection{Crowd dynamics: fast exit scenario}\label{sec:pedexit}
In this section we consider a model for crowd dynamics taken
from Reference~\cite{BdFMW13}.
We set the model in the spatial domain $\Omega=[-1,1]$,
whose boundaries represent the exit doors. The
non-local interaction kernel $p(x,y)$ is null and the
selective function $s(t,x,\rho)$
is $1-\rho$ (hence, we employ the exponential Euler method~\eqref{eq:expEuler}
for the forward equation).
The diffusion parameter is $\sigma=\sqrt{0.04}$ and
the exit intensity flux is $\beta=10$.
The initial density function
models the presence of two distinct groups, namely
$\rho_0(x) = 0.9\ee^{-100(x+0.4)^2}+0.65\ee^{-150x^2}$.

Similarly to the opinion dynamics case,
we first show that the expected temporal rate of convergence
of the exponential integrators is preserved after a complete
solution of the model. To this purpose,
we discretize this problem with $n=200$ spatial discretization points
and with different number of time steps, from $m=300$ to $m=700$,
up to the final time $T=2$.
After the stabilization of the functional $\mathcal{J}$ in the steepest
descent algorithm,
we measure the error at final time for 
$\brho(t)$ and at initial time for
$\bpsi(t)$ with respect to reference solutions.
We display in Figure~\ref{fig:orderpedexit} the obtained relative errors
which again confirm the expected accuracy and rate of convergence.

\begin{figure}[htb!]
  \centering
%
%
\begin{tikzpicture}

\begin{axis}[%
width=1.464in,
height=1.464in,
at={(0.769in,0.477in)},
scale only axis,
xmode=log,
xmin=300,
xmax=700,
xminorticks=true,
xtick={300,400,500,600,700},
xticklabels={{300},{400},{500},{600},{700}},
xlabel style={font=\color{white!15!black}},
xlabel={$m$},
ylabel style={font=\color{white!15!black}},
ylabel={Relative error},
ymode=log,
ymin=1e-3,
ymax=1e-2,
yminorticks=true,
axis background/.style={fill=white},
]
\addplot [color=red, line width=1.5pt, only marks, mark size=4pt, mark=x, mark options={solid, red}, forget plot]
  table[row sep=crcr]{%
300	0.00311925860975009\\
400	0.00231084992834468\\
500	0.00182194508005949\\
600	0.00149443989054837\\
700	0.00125975292908461\\
};
\addplot [color=black, dashed, line width=1.1pt, forget plot]
  table[row sep=crcr]{%
300	0.00293942350119743\\
700	0.00125975292908461\\
};
\end{axis}

\begin{axis}[%
width=1.464in,
height=1.464in,
at={(3.075in,0.477in)},
scale only axis,
xmode=log,
xmin=300,
xmax=700,
xminorticks=true,
xtick={300,400,500,600,700},
xticklabels={{300},{400},{500},{600},{700}},
xminorticks=true,
xlabel style={font=\color{white!15!black}},
xlabel={$m$},
ylabel style={font=\color{white!15!black}},
ylabel={Relative error},
ymode=log,
ymin=1e-3,
ymax=1e-2,
yminorticks=true,
axis background/.style={fill=white},
]
\addplot [color=blue, line width=1.5pt, only marks, mark size=4pt, mark=o, mark options={solid, blue}, forget plot]
  table[row sep=crcr]{%
300	0.00678057802270684\\
400	0.00501861461886648\\
500	0.00395655685095689\\
600	0.00324742845398448\\
700	0.00274066980165971\\
};
\addplot [color=black, dashed, line width=1.1pt, forget plot]
  table[row sep=crcr]{%
300	0.00639489620387265\\
700	0.00274066980165971\\
};
\end{axis}

\end{tikzpicture}%
  \caption{Relative errors in infinity norm
    of $\brho(T)$ (left, $T=2$) and $\bpsi(0)$ (right),
    with respect to a reference solution, for the
    pedestrian model described in Section~\ref{sec:pedexit}
    with $n=200$
    spatial discretization points and varying number of time steps $m$.
    The reference line of order 1 is also displayed.}
\label{fig:orderpedexit}
\end{figure}

\begin{figure}[htb!]
  \centering
  \includegraphics[scale=0.22]{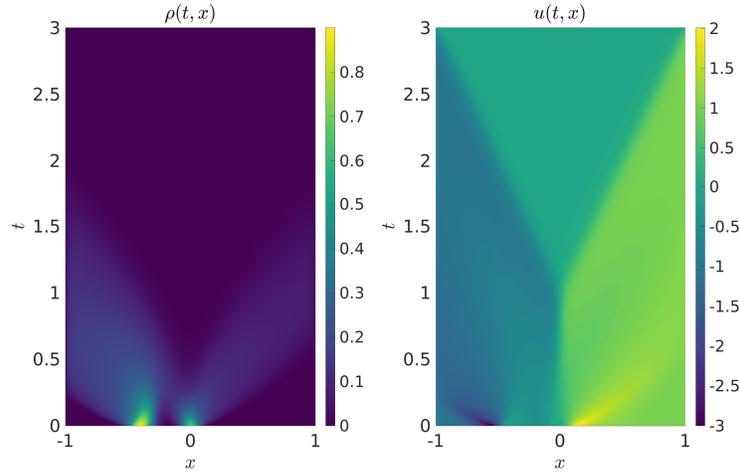}
  \caption{Evolution of the density $\rho(t,x)$ (left) and of the
    control $u(t,x)$ (right) up to final time $T=3$ for the
    two-group crowd model
    described in Section~\ref{sec:pedexit} with $n=1000$ spatial
    discretization points and $m=250$ time steps.}
\label{fig:behaviorpedexit}
\end{figure}

\begin{figure}[htb!]
  \centering
%
%
%
\begin{tikzpicture}

\begin{axis}[%
width=3.585in,
height=1.932in,
at={(0.769in,0.477in)},
scale only axis,
xmin=0,
xmax=14,
xlabel style={font=\color{white!15!black}},
xlabel={$\ell$},
ymin=0.22,
ymax=0.38,
ylabel style={font=\color{white!15!black}},
ylabel={$\mathcal{J}(u^\ell)$},
axis background/.style={fill=white},
legend style={legend cell align=left, align=left, draw=white!15!black}
]
\addplot [color=blue, only marks, line width=1.5pt, mark size=4pt, mark=x, mark options={solid, blue}]
  table[row sep=crcr]{%
1	0.373513241196724\\
2	0.372801263801694\\
3	0.367440579404244\\
4	0.356018699007276\\
5	0.336501997508616\\
6	0.312434843194177\\
7	0.290937394007775\\
8	0.272145916854953\\
9	0.25495668810539\\
10	0.244094186853728\\
11	0.23785047220661\\
12	0.233569518883002\\
13	0.230757325165227\\
14	0.229317747595152\\
};
\addlegendentry{Functional value}

\addplot [color=red, line width=1.5pt]
  table[row sep=crcr]{%
0	0.229317747595152\\
14	0.229317747595152\\
};
\addlegendentry{Reached value: 0.2293}

\end{axis}
\end{tikzpicture}%
  \caption{Value of the functional $\mathcal{J}(u^\ell)$ at the successive
iterations of the steepest descent method for the two-group crowd model
described in Section~\ref{sec:pedexit} ($n=1000$ and $m=250$).}
\label{fig:functionalpedexit}
\end{figure}
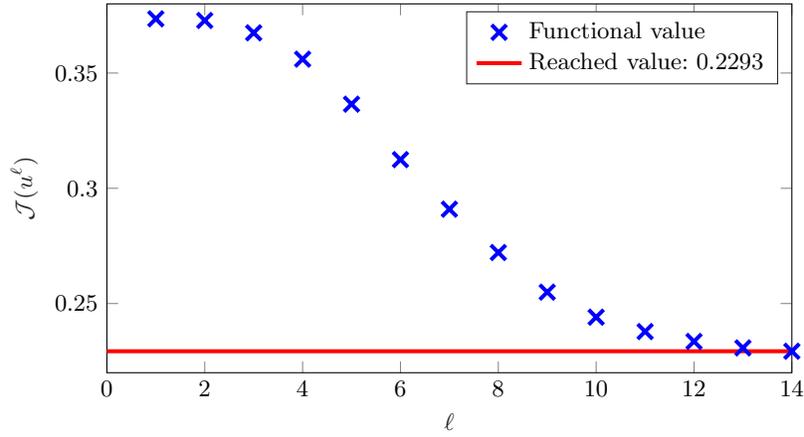

Then, we solve the same model up to the final time $T=3$ and
show its behavior.
We discretize this problem with $n=1000$ spatial discretization points
and $m=250$ time steps.
We show the evolution of the density  and of the control in
Figure~\ref{fig:behaviorpedexit}, where we can clearly see the exit of
the crowd from the two doors.
Moreover, we show in Figure~\ref{fig:functionalpedexit} the value
of the functional $\mathcal{J}(u^\ell)$ at the successive
iterations of the steepest descent method. We observe that
the method needs 14 iterations to reach the input tolerance $2\cdot 10^{-3}$.
Finally, the overall computational time of this simulation is
about 45 seconds.

\subsection{Mass transfer problem via optimal control}\label{sec:distmatch}
\begin{figure}[htb!]
  \centering
  \input{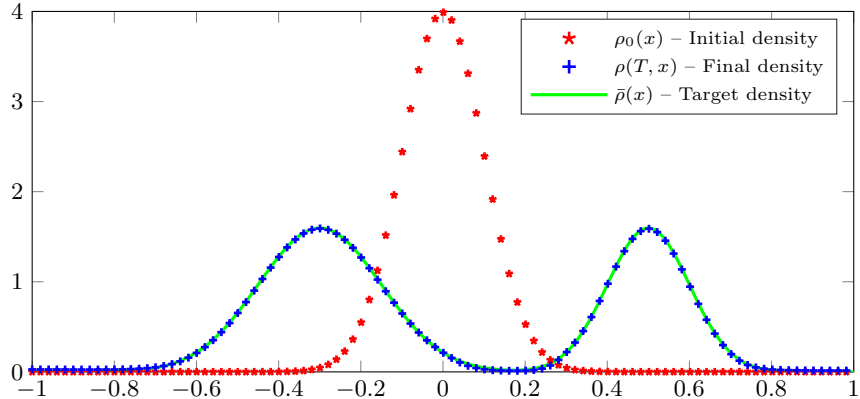}\\
  \caption{Density functions at initial time and at final time
    for the mass transfer problem
    described in Section~\ref{sec:distmatch} with $n=1000$ spatial
    discretization points and $m=200$ time steps.}
\label{fig:behaviordistmatch}
\end{figure}

In this final example, we present an optimal control approach to a mass transfer
problem, see for instance References~\cite{benamou2000computational,santambrogio2015optimal},
where the particle density accounts for non-local
interactions~\cite{bongini2017optimal,carrillo2012confinement}.
Hence, the goal is to move the initial density function in the spatial domain $\Omega=[-1,1]$
\begin{equation*}
  \rho_0(x) = C(\ee^{-(x-\mu_0)^2/(2\sigma_0^2)}),
\end{equation*}
where $\mu_0=0$, $\sigma_0=0.1$, and $C$ is defined so that
$\int_\Omega\rho_0(x)dx=1$,
to a target one
\begin{equation*}
  \bar{\rho}(x)={\bar C}\left(\ee^{-(x-\mu_1)^2/(2\sigma_1^2)}+
  \ee^{-(x-\mu_2)^2/(2\sigma_2^2)}\right),
\end{equation*}
where $\mu_1=0.5$, $\sigma_1=0.1$, $\mu_2=-0.3$, and $\sigma_2=0.15$,
and $\bar C$ is defined so that $\int_\Omega\bar\rho(x)dx=1$.
The boundary conditions are of zero-flux type, the running cost
is $e(t,x,\rho) = |\rho-\bar{\rho}|^2$,
the interaction kernel is of Sznajd type
$p(x,y)=(x^2-1)/20$, and the selective function is $s(t,x,\rho)=1$.
The penalization parameter is $\gamma=0.1$ and the diffusion parameter
is $\sigma=\sqrt{0.02}$.
We discretize the problem with $n=1000$ spatial grid points and
$m=200$ time steps, and we run the simulation up to the final time $T=3$.
We consider  a terminal cost given by
$c(T,x,\rho(T,x))=\lvert\rho(T,x)-\bar{\rho}(x)\rvert^2$, which translates
into
$\psi_T(x)=\rho(T,x)-\bar{\rho}(x)$. In
Figure~\ref{fig:behaviordistmatch} we plot the density functions
at the initial and the final time, and we can observe that the 
initial density is correctly transported to the target one. In addition, in
Figure~\ref{fig:rhoudistmatch} we present the evolution of the density and
of the control. 
Finally, we show in Figure~\ref{fig:functionaldistmatch} the values
of the functional $\mathcal{J}(u^\ell)$ at the successive
iterations of the steepest descent method. We observe that
the method needs 33 iterations to reach the input tolerance $2\cdot 10^{-3}$,
with an overall computational time of this simulation of roughly 75 seconds.

\begin{figure}[htb!]
  \centering
  \includegraphics[scale=0.22]{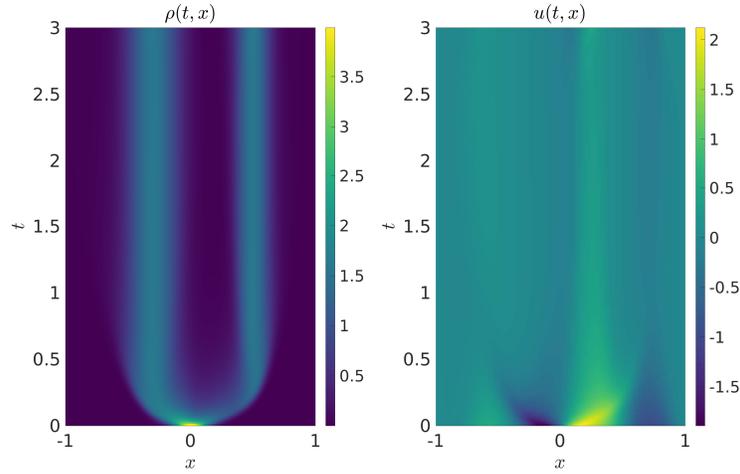}
  \caption{Evolution of the density $\rho(t,x)$ (left) and of the
    control $u(t,x)$ (right) up to final time $T=3$ for the
    mass transfer problem
    described in Section~\ref{sec:distmatch} with $n=1000$ spatial
    discretization points and $m=200$ time steps.}
\label{fig:rhoudistmatch}
\end{figure}
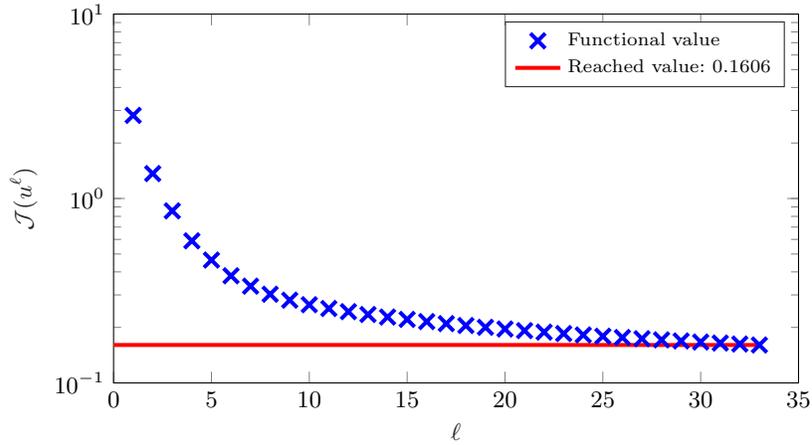
\begin{figure}[htb!]
  \centering
%
%
%
\begin{tikzpicture}

\begin{axis}[%
width=3.585in,
height=1.932in,
at={(0.758in,0.481in)},
scale only axis,
xmin=0,
xmax=35,
xlabel style={font=\color{white!15!black}},
xlabel={$\ell$},
ymin=0.1,
ymax=10,
ymode=log,
ylabel style={font=\color{white!15!black}},
ylabel={$\mathcal{J}(u^\ell)$},
axis background/.style={fill=white},
legend style={legend cell align=left, align=left, draw=white!15!black, font=\scriptsize}
]
\addplot [color=blue, line width=1.5pt, only marks, mark=x, mark size=4pt, mark options={solid, blue}]
  table[row sep=crcr]{%
1	2.8232200495847\\
2	1.36499697240198\\
3	0.856832276199011\\
4	0.589411490398348\\
5	0.463747782910922\\
6	0.381020840961113\\
7	0.334842541706958\\
8	0.302612684123614\\
9	0.28107589272649\\
10	0.265500261613218\\
11	0.253282539026929\\
12	0.243464484263135\\
13	0.234966819129323\\
14	0.227610021315361\\
15	0.221006078146646\\
16	0.215085729070666\\
17	0.209691043302795\\
18	0.204768384022931\\
19	0.200238952116048\\
20	0.196059658535403\\
21	0.192183865735593\\
22	0.188578697213063\\
23	0.185213217217325\\
24	0.182062798574176\\
25	0.179105435484739\\
26	0.176322590576763\\
27	0.173697916389736\\
28	0.171217256184661\\
29	0.168868114622628\\
30	0.166639494406984\\
31	0.164521635267856\\
32	0.162505854916082\\
33	0.160584397273363\\
};
\addlegendentry{Functional value}

\addplot [color=red, line width=1.5pt]
  table[row sep=crcr]{%
0	0.160584397273363\\
33	0.160584397273363\\
};
\addlegendentry{Reached value: 0.1606}

\end{axis}
\end{tikzpicture}%
  \caption{Value of the functional $\mathcal{J}(u^\ell)$ at the successive
    iterations of the steepest descent method for the mass transfer problem
    described in Section~\ref{sec:distmatch}.}
\label{fig:functionaldistmatch}
\end{figure}
\section{Conclusions}\label{sec:conclusions}
We presented a mean-field optimal control model where the constraint is
represented by a nonlinear PDE with non-local interaction term and
diffusion describing the evolution of a continuum of agents.
We provide, at a formal level, first order optimality conditions, resulting in
a forward-backward coupled system with associated boundary conditions.
Thus, a reduced gradient method is derived for the synthesis of the
mean-field control, where the primal and adjoint equations are efficiently
solved by using exponential integrators. Our proposed approach has been
successfully tested on various examples from the literature,
including models of opinion formation and pedestrian dynamics in
the one-dimensional setting.
In future works we plan to exploit the efficiency of exponential integrators
to tackle higher dimensional problems (possibly using 
ad hoc techniques for tensor structured
problems~\cite{CCEOZ22,CCZ22c,CCZ22b})
and scenarios where a fine 
spatial discretization is required to correctly capture the behavior of
the controlled dynamics.
\subsubsection*{Acknowledgments} The authors were partially supported by the MIUR-PRIN Project 2017,
No. 2017KKJP4X \emph{Innovative numerical methods for evolutionary partial differential equations
and applications}, and by RIBA 2019, No. RBVR199YFL \emph{Geometric Evolution of Multi Agent
Systems}. 
\bibliographystyle{spmpsci}
\bibliography{biblio_mfexpint.bib}
\end{document}